\documentclass[10pt,reqno]{amsart}

\usepackage{amssymb, amsmath}
\usepackage{hyperref}
\newcommand{\1}{\textbf{1}}
\renewcommand{\a}{\alpha}
\newcommand{\A}{\textit{\textbf{u}}}
\renewcommand{\b}{\beta}
\newcommand{\B}{\textit{\textbf{v}}}
\renewcommand{\c}{\gamma}

\newcommand{\I}{\tilde{I}}
\newcommand{\J}{\tilde{J}}

\newcommand{\m}{\mathfrak{m}}
\renewcommand{\o}{\otimes}

\newcommand{\pb}{\textrm{\textbf{p}}}
\newcommand{\s}{\sigma}
\renewcommand{\sb}{\textrm{\textbf{s}}}
\renewcommand{\ss}{\overline{\sigma}}
\renewcommand{\t}{\tau}
\renewcommand{\tt}{\overline{\tau}}

\newcommand{\x}{\textit{\textbf{x}}}

\newcommand{\y}{\textit{\textbf{y}}}
\newcommand{\z}{\textit{\textbf{z}}}
\newcommand{\bb}[1]{\mathbb{#1}}
\renewcommand{\rm}[1]{\textrm{#1}}
\newcommand{\mc}[1]{\mathcal{#1}}

\newtheorem{theorem}{Theorem}[section]

\newtheorem*{definition*}{Definition}

\newtheorem{lemma}[theorem]{Lemma}
\newtheorem{corollary}[theorem]{Corollary}
\newtheorem*{toreq}{Theorem \ref{torequiv}}
\newtheorem*{amit}{Theorem \ref{thmamit}}
\newtheorem*{quot}{Corollary \ref{torquot}}

\theoremstyle{remark}

\newtheorem{example}[theorem]{Example}

\begin{document}

\title{Affine Toric Equivalence Relations are Effective}
\author{Claudiu Raicu}
\date{\today}
\address{Department of Mathematics, University of California, Berkeley\newline
\indent Institute of Mathematics "Simion Stoilow" of the Romanian Academy}
\email{claudiu@math.berkeley.edu}

\begin{abstract} Any map of schemes $X\to Y$ defines an equivalence relation $R=X\times_Y X\to X\times X$, the relation of ``being in the same fiber". We have shown elsewhere that not every equivalence relation has this form, even if it is assumed to be finite. By contrast, we prove here that every toric equivalence relation on an affine toric variety does come from a morphism and that quotients by finite toric equivalence relations always exist in the affine case. In special cases, this result is a consequence of the vanishing of the first cohomology group in the Amitsur complex associated to a toric map of toric algebras. We prove more generally the exactness of the Amitsur complex for maps of commutative monoid rings.
\end{abstract}
\maketitle

\section{Introduction}

The question that motivated this paper was the following, asked by J\'anos Koll\'ar: given an $S$-scheme $X$ and a finite scheme theoretic equivalence relation $R\subset X\times_S X$, does there exist an $S$-scheme $Y$ and a finite surjective map $X\to Y$ over $S$ such that $R\simeq X\times_Y X$? The answer is in general negative and we gave an example of this phenomenon in an appendix to \cite{Kol}; it is reproduced as Example \ref{exnoneff} below.

Recall (\cite{FGAex}, \cite{Kol}) that given a scheme $X$ over a base $S$, a \textit{scheme theoretic equivalence relation} on $X$ over $S$ is an $S$-scheme $R$ together with a morphism $f:R\to X\times_S X$ over $S$ such that for any $S$-scheme $T$, the set map $f(T):R(T)\to X(T)\times X(T)$ is injective and its image is the graph of an equivalence relation on $X(T)$ (where for $S$-schemes $Z,T$ we denote by $Z(T)$ the set of $S$-maps from $T$ to $Z$). An equivalence relation is \textit{finite} if the two maps $f_1,f_2:R\rightrightarrows X$ corresponding to $f$ are finite, in which case $f$ is a closed immersion (see \cite{Kol}). In this paper we will only be interested in equivalence relations defined by closed subschemes of $X\times_S X$, so we may as well take $f$ to be a closed immersion in the previous definition. We say that an equivalence relation $R$ is \textit{effective} if there exists a morphism $X\to Y$ such that the induced map $R\to X\times_Y X$ is an isomorphism. We make the following

\begin{definition*} Let $k$ be a field and $X/k$ a (not necessarily normal) toric variety. A scheme theoretic equivalence relation on $X$ defined by a closed subscheme $R$ of $X\times_k X$ is called toric if $R$ is invariant under the diagonal action of the torus on $X\times_k X$.
\end{definition*}

With this definition, we can now state the main result of this paper.

\begin{toreq} Let $k$ be a field, $X/k$ an affine toric variety, and $R$ a toric equivalence relation on $X$. Then there exists an affine toric variety $Y$ together with a toric map $X\to Y$ such that $R\simeq X\times_Y X$.
\end{toreq}

\noindent This result cannot be generalized to arbitrary toric varieties, as shown in Example \ref{exnonaff}. Theorem \ref{torequiv} turns out to be related to the vanishing of the first cohomology group in the Amitsur complex (defined in Section \ref{amitsur})
\[C(A,B):\ B\to B\o_A B\to\cdots\to B^{\o_A m}\to\cdots\]
associated to a toric map of toric algebras $A\to B$. We prove the following more general

\begin{amit} Let $k$ be any commutative ring, $\t,\s$ commutative monoids, and $\varphi:\t\to\s$ a map of monoids. If $A=k[\t]$, $B=k[\s]$, and $B$ is considered as an $A$-algebra via the map $A\to B$ induced by $\varphi$, then the Amitsur complex $C(A,B)$ is exact.
\end{amit}

Questions about equivalence relations usually arise from the desire to construct a good quotient scheme. A morphism $q:X\to Q$ is called the \textit{categorical quotient} of $X$ by the equivalence relation $R$ if it is a coequalizer in the category of $S$-schemes of the morphisms $f_1,f_2:R\rightrightarrows X$. It is \textit{effective} if the induced map $R\to X\times_Q X$ is an isomorphism. The categorical quotient $q$ is a \textit{geometric quotient} if it is finite and for every geometric point $\rm{Spec }K\to S$, the fibers of $q_K:X_K(K)\to Q_K(K)$ are the $f(R_K(K))$-equivalence classes of $X_K(K)$. A finite effective categorical quotient is also a geometric quotient (see \cite{Kol}).

The question of existence of quotients $X/R$ when $R$ is a scheme (or set) theoretic equivalence relation on $X$ is complicated in general. In the proper flat case, the existence of categorical quotients is the result of work by Grothendieck (\cite{FGA}) and Altman and Kleiman (\cite{AlKl}). For a more recent account of this result see \cite{FGAex}. Finite set theoretic equivalence relations were studied by Koll\'ar, who proved the existence of geometric quotients in the case of schemes of finite type over a field of positive characteristic (\cite{Kol}). In the affine case, Olivier gave a necessary and sufficient condition for the existence of categorical quotients in \cite{Oli}, which we recall as Lemma \ref{olivier}.

Unfortunately, the categorical quotient need not exist for an arbitrary affine toric equivalence relation, as Example \ref{exnonfinite} shows. However, if the equivalence relation $R$ on $X$ is both finite and toric, then combining Theorem \ref{torequiv} with Lemma \ref{olivier} we get

\begin{quot} Let $k$ be a field, $X/k$ an affine toric variety, and $R$ a finite toric equivalence relation on $X$. Then the geometric quotient $X/R$ exists and is effective.
\end{quot}

The structure of this paper is as follows. In Section \ref{prelim} we translate the definitions of equivalence relations in the affine case, and introduce the basic notions about monoids that will be used throughout the rest of the paper. We then prove the exactness of the Amitsur complex for maps of monoid rings in Section \ref{amitsur}. The proof of the effectiveness of toric equivalence relations is the content of Section \ref{atoric}. Finally, Section \ref{quots} deals with the existence of effective geometric quotients for finite toric equivalence relations.

\section{Preliminaries}\label{prelim}

\subsection{Equivalence relations in the affine case}\label{affequiv} Let $k$ be a field, $X=\bb{A}_k^n$ the $n$-dimensional affine space over $k$. Then $\mc{O}_X\simeq k[\x]$, where $\x=(x_1,\cdots,x_n)$. To give an equivalence relation $R\subset X\times_k X$ is equivalent to giving an ideal $I(\x,\y)\subset k[\x,\y]$ that satisfies the following properties:
\begin{enumerate}
  \item[i)] (reflexivity) $I(\x,\y)\subset (x_1-y_1,\cdots,x_n-y_n)$.
  \item[ii)] (symmetry) $I(\x,\y)=I(\y,\x).$
  \item[iii)] (transitivity) $I(\x,\z)\subset I(\x,\y)+I(\y,\z)$ in $k[\x,\y,\z]$.
\end{enumerate}
$R$ is finite if and only if $I$ satisfies
\begin{enumerate}
  \item[iv)] (finiteness) $k[\x,\y]/I(\x,\y)$ is finite over $k[\x]$.
\end{enumerate}

Analogue conditions hold if we replace $k[\x]$ with $k[\s]$ for a finitely generated submonoid $\s$ of $\bb{Z}^n$. The equivalence relation $R$ comes from a map to an affine $k$-variety $Y$ if and only if $I(\x,\y)$ is generated by differences $f_i(\x)-f_i(\y)$, $i=1,\cdots,m$, where the $f_i$'s are generators of the $k$-algebra $\mc{O}_Y$. In the appendix to \cite{Kol} we gave the following example of a finite noneffective equivalence relation.

\begin{example}\label{exnoneff} Let $k$ be any field, consider $X=\bb{A}_k^2$ and let $R$ be the equivalence relation on $X$ defined by the ideal
\[I(\x,\y)=(x_1^2-y_1^2,x_1x_2-x_2^2-y_1y_2+y_2^2,x_2^3-y_2^3,(x_1y_2-x_2y_1)y_2^3)\subset k[\x,\y].\]
Then $R$ is a finite noneffective equivalence relation. This phenomenon doesn't occur in the toric case, as will be explained in Section \ref{atoric}.
\end{example}

\subsection{Monoids}

By a \textit{monoid} we will understand a commutative semigroup with identity element $1$. We will write the monoid operation multiplicatively whenever we refer to general monoids (Section \ref{amitsur}), and additively when we deal with affine monoids (Section \ref{atoric}). Given monoids $\t,\s_1,\s_2$ with $\t$ mapping to both $\s_1$ and $\s_2$ we can form the tensor product (see \cite{How}, Section 8.1)
\[\s_1\o_{\t}\s_2:=\s_1\times\s_2/\sim,\]
where $\sim$ is the equivalence relation on $\s_1\times\s_2$ generated by the relation
\[R=\{((s_1t,s_2),(s_1,ts_2)):s_1\in\s_1,s_2\in\s_2,t\in\t\}.\]

Given a commutative ring $k$, we can form the monoid rings $T=k[\t]$, $S_i=k[\s_i]$ and their tensor product $S_1\o_T S_2$. We let $S_{12}=k[\s_1\o_{\t}\s_2]$. The natural maps $\s_1,\s_2\to\s_1\o_{\t}\s_2$ induce maps $S_1,S_2\to S_{12}$ which agree when restricted to $T$. We obtain a natural isomorphism $S_1\o_T S_2\simeq S_{12}$ and make the convention to identify the two constructions via this isomorphism throughout the paper.

Given a map of monoids $\t\to\s$ with associated ring map $T\to S$ we can form the $n$-fold tensor products $\s^{\o_{\t} n}=\s\o_{\t}\s\o_{\t}\cdots\o_{\t}\s$ and $S^{\o_{T} n}=k[\s^{\o_{\t} n}]$. We write a general element of $\s^{\o_{\t} n}$ as $\sb=s_1\o_{\t} s_2\o_{\t}\cdots\o_{\t} s_n$ and call it a \textit{monomial}. There is a natural surjection $\s^n=\s^{\o_{\langle 1\rangle}n}\twoheadrightarrow\s^{\o_{\t} n}$ which we will use to choose particular representatives for elements of $\s^{\o_{\t} n}$.

We define multiplication maps $\mu_{j}^{i}:\s^{\o_{\t} n}\to\s^{\o_{\t}(n-1)}$ by
\[\mu_{j}^i(\sb)=s_1\o_{\t}\cdots\o_{\t} s_{i-1}\o_{\t} s_is_j\o_{\t} s_{i+1}\o_{\t}\cdots\widehat{s_j}\cdots\o_{\t} s_n\]
for $i\neq j\leq n$. They are well-defined because $\t$ and $\s$ are commutative. We also define degeneracy maps $\xi_i:\s^{\o_{\t} n}\to\s^{\o_{\t}(n+1)}$ by
\[\xi_i(\sb)=s_1\o_{\t}\cdots\o_{\t} s_{i-1}\o_{\t} 1\o_{\t} s_i\o_{\t}\cdots\o_{\t} s_n\]
for $i=1,2,\cdots,n+1$. We write $\mu_{j}^i,\xi_i$ also for the corresponding induced maps on $S^{\o_T n}$.

The \textit{degree} of a monomial $\sb=s_1\o_{\t}\cdots\o_{\t} s_n\in\s^{\o_{\t} n}$ is defined by
\[\rm{deg}(\sb)=s_1\cdots s_n\in\s.\]
This is the same as $(\mu_{2}^1)^{n-1}(\sb)$, and hence does not depend on the particular representation of the monomial $\sb$.

In general, we will omit to write $\t$ whenever it's understood. Also, we will make no notational distinction between an element of $\s^{\o_{\t} n}$ and its representative in $\s^n$, as this will usually be irrelevant and/or clear from the context.

\paragraph{\textit{Affine monoids}}

\begin{definition*}
An affine monoid is a finitely generated cancellative torsion free monoid, or equivalently, a finitely generated submonoid of a finitely generated free abelian group.
\end{definition*}

For such monoids, we write the monoid operation additively and denote by 0 the identity element. Given an affine monoid $\s$ we write $k[\x]$ for the monoid ring $k[\s]$, and $\x^s$ for the monomial corresponding to $s\in\s$. Similarly we write $k[\x,\y]$, $k[\x,\y,\z]$ for $k[\s^2]$, $k[\s^3]$ (these will be the only powers of $\s$ we will be interested in).

We let $\rm{grp}(\s)$ be the group generated by $\s$, and $l(\s)$ the largest subgroup of $\s$. We say that $\s$ is \textit{pointed} if $l(\s)=0$. A pointed affine monoid admits a \textit{monomial order}, that is a well ordering of its elements which is compatible with addition in the sense that $u+w\leq v+w$ whenever $u\leq v$. Given an arbitrary affine monoid $\s$, we consider the equivalence relation on $\s$ given by $u\sim v$ if and only if $u-v\in l(\s)$ (with the difference computed inside $\rm{grp}(\s)$), and let $\overline{\s}=\s/\sim$. We define the \textit{reduced degree} (\textit{$r$-degree} for short) of a monomial $\sb\in\s^{\o_{\t} n}$ (where $\t$ is a monoid mapping to $\s$) to be the class of $\rm{deg}(\sb)$ in $\overline{\s}$.

Notice that $\overline{\s}$ is itself a cancellative monoid without invertible elements, and although it might not be torsion free ($\overline{\s}$ being torsion free is equivalent to $l(\s)$ being a direct summand in $\rm{grp}(\s)$) it admits a partial order compatible with addition as in the following

\begin{lemma}\label{lmords} With the above notations, $\overline{\s}$ admits a partial order for which any decreasing chain stabilizes and with the property that if $u,v,w\in\overline{\s}$ are such that $u+v=w$, then $u,v\leq w$, and the inequality is strict unless $u$ or $v$ is zero.
\end{lemma}

\begin{proof} Let $l$ be the direct summand of $\rm{grp}(\s)$ containing $l(\s)$ as a subgroup of finite index. Consider the equivalence relation on $\s$ given by $u\sim' v$ if and only if $u-v\in l$, and let $\widetilde{\s}$ denote the quotient $\s/\sim'$. $\widetilde{\s}$ is a pointed affine monoid, hence it admits a monomial order. Let $\pi$ denote the natural projection $\overline{\s}\to\widetilde{\s}$ and define the partial order on $\overline{\s}$ via $u>v$ whenever $\pi(u)>\pi(v)$. The descending chain condition on $\overline{\s}$ and the fact that $u+v=w$ can only happen when $u,v\leq w$ follow from the corresponding properties of $\widetilde{\s}$ and the fact that $\pi(u)=0\in\widetilde{\s}$ if and only if $u=0\in\overline{\s}$.

\end{proof}

The existence of an order as in the preceding lemma will be important in the inductive argument of Section \ref{toricgenlcase}, and will allow us to apply Nakayama's lemma in the proof of Corollary \ref{torquot}.

\section{Exactness of the Amitsur complex for maps of monoid rings}\label{amitsur}

For a commutative ring $A$ and an $A$-algebra $B$, we consider the Amitsur complex
\[C(A,B):\ B\to B\o_A B\to\cdots\to B^{\o_A m}\to\cdots\]
starting in degree zero, and with differentials given by the formula
\[d_{m-1}(b_1\o b_2\o\cdots\o b_m)=\sum_{i=1}^{m+1}(-1)^i b_1\o\cdots\o b_{i-1}\o 1\o b_i\o\cdots\o b_m.\]
It is well known that if $B$ is a faithfully flat or augmented $A$-algebra, then $C(A,B)$ is exact (see \cite{Jong}, Lemma 2.6). The main result in this section is that $C(A,B)$ is also exact for maps of monoid rings. Notice that Example \ref{exnoneff} shows that exactness of the Amitsur complex fails for arbitrary ring extensions.

\begin{theorem}\label{thmamit} Let $k$ be any commutative ring, $\t,\s$ commutative monoids, and $\varphi:\t\to\s$ a map of monoids. If $A=k[\t]$, $B=k[\s]$, and $B$ is considered as an $A$-algebra via the map $A\to B$ induced by $\varphi$, then the Amitsur complex $C(A,B)$ is exact.
\end{theorem}

An important difference between the monoidal case and the faithfully flat and augmented cases is that in the latter cases $\rm{Ker }d_0=A$, while in the former case $\rm{Ker }d_0$ is usually bigger than $A$. Though we will not use this here, Isbell's Zigzag Theorem gives a nice criterion for when the equality $\rm{Ker }d_0=A$ holds in the monoidal case (see \cite{HowIs}, or \cite{How}, Theorem 8.3.4).

We fix $\s,\t$ as in the statement of Theorem \ref{thmamit}. We give some preliminary results and definitions before proceeding to the proof of the theorem.

\begin{definition*} A monomial $\sb=s_1\o\cdots\o s_n\in\s^n$ is said to be normalized if it has the maximal number of $s_i$'s equal to 1 among all representatives of its class in $\s^{\o_{\t}n}$.
\end{definition*}

\begin{lemma}\label{lmnorm} A monomial $\sb\in\s^n$ is normalized if and only if $\xi_i(\sb)$ is normalized for some (all) $i\in\{1,\cdots,n+1\}$.
\end{lemma}

\begin{proof} The `if' part is clear. To prove the converse, notice that $\sb=\mu_{i+1}^i(\xi_i(\sb))$ (or $\mu_{n+1}^n(\xi_{n+1}(\sb))$ if $i=n+1$), and $\xi_i(\sb)$ has precisely one more 1 than $\sb$. Any multiplication map can decrease the number of 1's in a monomial by at most one, hence the conclusion follows.

\end{proof}

\begin{lemma}\label{lm1s} If $\sb,\sb'$ are normalized and equal in $\s^{\o_{\t} n}$, then either each of them has only one factor different from 1 (which must then be the same for both, equal to the degree), or they have the factors equal to 1 in the same positions.
\end{lemma}

\begin{proof} Suppose there exists some $i\in\{1,\cdots,n\}$ for which one and only one of $s_i,s_i'$ is equal to 1, say $i=1$ and $s_1=1$. Suppose then that $s_j'\neq 1$ for some $j>1$. Then since $\sb=\xi_1(\mu_{1}^j(\sb))$ the same must be true for $\sb'$, but $\xi_1(\mu_{1}^j(\sb'))$ has one more $1$ than $\sb'$, contradicting the fact that $\sb'$ was normalized. Therefore $\sb'$ has at most one factor different from $1$ and the same holds for $\sb$ by applying a similar argument with $s_1$ replaced by whichever $s_i'=1$ has the property that $s_i\neq 1$. Since $\sb,\sb'$ are equal, in particular they have the same degree, i.e. $s_1\cdots s_n=s_1'\cdots s_n'=s$ in $\s$, hence both $\sb,\sb'$ are products involving precisely one $s$ and $(n-1)$ 1's.

\end{proof}

We partially order the elements of $\s^n$ by saying that $\sb>\sb'$ if there exists some $i$ for which $s_i,s_i'$ are not both equal to, or both different from 1, and the first time this happens we have $s_i\neq 1$. Given an element $f=\sum a_{\sb}\sb$ we call \textit{dominant} any term $a_{\sb}\sb$ with $a_{\sb}\neq 0$ and $\sb$ maximal with respect to this order. We call the \textit{type} of a monomial $\sb$ the set of indices $i$ for which $s_i\neq 1$, and order the possible types accordingly.

\begin{proof}[Proof of Theorem \ref{thmamit}]

Let $f=\sum a_{\sb}\sb\in k[\s^{\o_{\t} n}]$ be an $(n-1)$-cocycle in $C(A,B)$, $n\geq 2$, and assume that all monomials in $f$ are normalized. We prove by double induction on the type and number of dominant terms that $f$ is a coboundary. Let $a_{\sb}\sb$ be a dominant term of $f$, and suppose that $\sb$ ends in an odd number of 1's, $\sb=s_1\o\cdots\o s_r\o 1\o\cdots\o 1$ with $n-r$ odd and $s_r\neq 1$. Let $\pb=\mu_{n}^{n-1}(\sb)$ and notice that $\sb>\xi_i(\pb)$ for $i=1,\cdots,r$, and $\sb=\xi_i(\pb)$ for $i>r$. Since $n-r$ is odd, it follows that $d(a_{\sb}\pb)=\pm a_{\sb}\sb+\{lower\ terms\}$, hence by induction $f\mp d(a_{\sb}\pb)$ is a coboundary, and the same is true for $f$. Notice that we are using Lemma \ref{lmnorm} to make sure that the terms in the new cocycle $f\mp d(a_{\sb}\pb)$ are normalized, or alternatively, to show that normalization can be done without increasing the type of a monomial.

Suppose now $\sb$ is as above, with $n-r$ even. It follows that $d(\sb)=\pm\sb\o 1+\{lower\ terms\}$, and since $df=0$, a combination of the following must hold
\begin{enumerate}
\item[i)] $\sb\o 1=\xi_i(\sb')$ for some $\sb'\neq\sb$ appearing in the expression of $f$.

\item[ii)] $\sb\o 1=\xi_i(\sb)$ for some $i\leq r$.
\end{enumerate}
If the latter case holds, then since $\sb\o 1$ and $\xi_i(\sb)$ are normalized (Lemma \ref{lmnorm}) and don't have the same type, it follows from Lemma \ref{lm1s} that $\sb=1\o\cdots\o s\o\cdots\o 1$ for some $s\in\s$. Applying $(\mu_3^2)^{n-r}\circ(\mu_2^1)^{r-1}$ to the equality $\sb\o 1=\xi_i(\sb)$, we get that $1\o s=s\o 1$. In particular all monomials involving one $s$ and $(n-2)$ 1's are equal, and we denote by $\pb$ their common value. If $n$ is odd, then $\sb=d(\pb)$ and we can conclude by induction as before, otherwise $d(\sb)=\sb\o 1$ and it suffices to deal with case i).

Suppose now i) holds. If $i=n$ or $n+1$ then applying $\mu_{n+1}^n$ to $\sb\o 1=\xi_{i}(\sb')$ we get $\sb=\sb'$, a contradiction, so we can assume $i\leq n-1$. By Lemma \ref{lmnorm}, $\sb\o 1$ and $\xi_i(\sb')$ are normalized. If they have only one factor $s\neq 1$, situated in distinct positions, it follows that $s\o 1=1\o s$ as before. But $\sb,\sb'$ are then monomials involving precisely one $s$ and $(n-1)$ 1's, and therefore equal, a contradiction. Otherwise, by Lemma \ref{lm1s} they must have the 1's in the same positions. Then $s_i=1$ and since $\sb\nless\sb'$, we must have $s_i'=1$. But $s_i'$ is the $(i+1)$-st factor of $\xi_i(\sb')$, hence $s_{i+1}=1$ and continuing in the same fashion we get $s_j=s_j'=1$ for all $j\geq i$. It follows that $\sb'=\mu_{n+1}^n(\xi_i(\sb'))=\mu_{n+1}^n(\sb\o 1)=\sb$, again a contradiction which concludes the proof of the theorem.

\end{proof}

\section{Affine toric equivalences}\label{atoric}

In this section $k$ is a field, $T$ a split algebraic torus over $k$, with character lattice $M\simeq\bb{Z}^n$, $\s$ a finitely generated submonoid of $M$ and $X=\rm{Spec }k[\s]$ the corresponding affine toric variety. We prove the following

\begin{theorem}\label{torequiv} Let $k$ be a field, $X/k$ an affine toric variety, and $R$ a toric equivalence relation on $X$. Then there exists an affine toric variety $Y$ together with a toric map $X\to Y$ such that $R\simeq X\times_Y X$.
\end{theorem}

We start by proving the theorem in the case when $\s$ is a pointed monoid, as an application of Theorem \ref{thmamit}. This gives a flavor of the general argument and provides a short proof of Theorem \ref{torequiv} in the case when $X=\bb{A}^n$. The general case is more complicated, and we consider it afterwards.

Unfortunately, we cannot expect to generalize Theorem \ref{torequiv} to an arbitrary toric variety $X$, as the following example shows.
\begin{example}\label{exnonaff}
Consider $X=\bb{P}^2$, $R=\Delta_X\cup(\bb{P}^1\times\bb{P}^1)$, where $\Delta_X$ denotes the diagonal and $\bb{P}^1$ is any torus-invariant line in $X$. If $R$ came from a map $f:X\to Y$ as above, then $f$ would have to contract the invariant $\bb{P}^1$, and therefore be constant.
\end{example}

We will write $k[\x]$ for $\mc{O}_X=k[\s]$ and $k[\x,\y]$, $k[\x,\y,\z]$ for $\mc{O}_{X\times_k X}$, $\mc{O}_{X\times_k X\times_k X}$ respectively. We denote the ideal of $R$ in $X\times_k X$ by $I=I(\x,\y)$. By assumption $I$ is homogeneous with respect to the $M$-grading induced by the diagonal action of $T$ on $X\times_k X$. Showing that the equivalence relation $R$ comes from a toric map to an affine toric variety is equivalent to proving that $I(\x,\y)$ is generated by binomial differences of the form $\x^w-\y^w$.

\subsection{The case where $\s$ is pointed}\label{spointed}

We fix a monomial ordering on $\s$ and a system of homogeneous generators of $I$. We will need the following

\begin{lemma}\label{lmcoc} Let $R$ be a toric equivalence relation on $X$ with ideal $I=I(\x,\y)$. If $f\in I$ is a homogeneous element and $\I$ is the ideal generated by the elements of $I$ of degree smaller than that of $f$, then $f$ satisfies the cocycle condition
\[f(\x,\y)+f(\y,\z)-f(\x,\z)\equiv 0\rm{ mod }\I(\x,\y)+\I(\y,\z)\subset k[\x,\y,\z].\]
\end{lemma}

\begin{proof} By transitivity of $R$ and homogeneity of $I$ we have that
\begin{equation}\label{eqncoc}
f(\x,\z)\equiv g(\x,\y)+h(\y,\z)\rm{ mod }\I(\x,\y)+\I(\y,\z)\subset k[\x,\y,\z],
\end{equation}
where $g,h\in I$ have the same degree as $f$. Letting $\y=\z$ and using reflexivity of $R$ we get that
\[f(\x,\z)\equiv g(\x,\z)\rm{ mod }\I(\x,\z)\subset k[\x,\z].\]
Similarly, letting $\x=\y$ we have
\[f(\y,\z)\equiv h(\y,\z)\rm{ mod }\I(\y,\z)\subset k[\y,\z].\]
This shows that we can replace $g$ and $h$ by $f$ in (\ref{eqncoc}) to get the desired conclusion.

\end{proof}

Assume now that $I$ is not generated by differences, and let $f\in I$ be a homogeneous generator of minimal degree which is not a difference. It follows that if $\I$ is the ideal generated by the elements of $I$ of degree smaller than that of $f$, then we must have $\I=(\x^{w_i}-\y^{w_i}:i=1,\cdots,r)$ for some $w_1,\cdots,w_r\in\s$. Letting $A=k[\x^{w_i}:i=1,\cdots,r]$, $B=k[\x]$ and using the previous lemma, we get that $f$ is a 1-cocycle in the Amitsur complex associated to the inclusion map $A\hookrightarrow B$. By Theorem \ref{thmamit}, $f$ is also a coboundary, i.e. it is congruent to a difference modulo $\I$. Since $f$ and $\I$ are homogeneous, it must be congruent to a binomial difference $c(\x^w-\y^w)$ modulo $\I$. If $c=0$ then $f\in\I$, otherwise we can replace it with $\x^w-\y^w$ in the system of generators of $I$ and conclude by induction.

\subsection{The general case}\label{toricgenlcase}

Given any ideal $J(\x,\y)\subset k[\x,\y]$ we will write $J(\x,\y,\z)$ for $J(\x,\y)+J(\y,\z)+J(\x,\z)\subset k[\x,\y,\z]$, and $J_0$ for the degree $0\in\s$ part of $J$. $k[\x,\y]_0$ is the coordinate ring of a torus $T'$ with character lattice $l(\s)$, diagonally embedded into $l(\s)^2$ via $s\mapsto(s,-s)$. We will write $\A^s$ for the character $\x^{s}\y^{-s}$ of $T'$ and similarly $\B^{s}$ for $\y^{s}\z^{-s}$. Elements $p(\x,\y)$ and $p(\x,\y,\z)$ of degree 0 will often be written as $p(\A)$ and $p(\A,\B)$ respectively.

Suppose that $I$ is not generated by differences and consider $\gamma\in\overline{\sigma}$ minimal with the property that $I$ is not generated by differences in $r$-degree at most $\gamma$ (we assume that $\overline{\s}$ is ordered as in Lemma \ref{lmords}). Let $\I$ be the ideal generated by the elements of $I$ of $r$-degree smaller than $\gamma$. It follows by the choice of $\gamma$ that $\I$ is generated by binomial differences $\x^{w_i}-\y^{w_i}$. We denote by $\t$ the submonoid of $\s$ generated by these $w_i$'s, so that $k[\x,\y]/\I(\x,\y)\simeq k[\s^{\o_{\t} 2}]$ and $k[\x,\y,\z]/\I(\x,\y,\z)\simeq k[\s^{\o_{\t} 3}]$. Let $(g_i)_{i=1,\cdots,m}$ generate $I$ modulo $\I$ in $r$-degree $\gamma$. We can assume that all $g_i$'s have the same degree, which we also denote by $\gamma$ (we will use the same notation for elements of $\s$ and $\overline{\s}$, as long as this doesn't cause confusion).

The strategy of proof is as follows:

\begin{itemize}

\item We show that we can assume $g_i=p_i\x^{\c}+q_i\y^{\c}$, for $p_i,q_i$ elements of degree $0\in\s$ in $k[\x,\y]$.

\item We prove that $J_0=(\I:\x^{\c})_0\subset k[\x,\y]_0=k[T']$ is the ideal of a subgroup scheme $T''$ of $T'$.

\item By studying the ideal generated by the $p_i$'s, we then reduce to the case when $p_1=1$ and $p_i=0$ for $i>1$.

\item In this case, we show that $\J=(q_2,\cdots,q_m)$ defines a subgroup scheme $T'''\subset T''$ and that $-q_1$ is a character of $T'''$. This is enough to conclude that $I$ is generated by differences in $r$-degree at most $\gamma$.

\end{itemize}

We start with the following

\begin{lemma}\label{lmcoc2} Suppose $g_i\in k[\x,\y]$, $i=1,\cdots,m$ are elements of degree $\gamma$ that generate $I/\I$ in $r$-degree $\gamma$. Then there exist elements $a_{ij},b_{ij}\in k[\x,\y,\z]$ of degree $0\in\s$ such that
\begin{equation}\label{coc2}
g_i(\x,\z)\equiv\sum_{j=1}^{m}(a_{ij}g_j(\x,\y)+b_{ij}g_j(\y,\z))\rm{ mod }\I(\x,\y,\z),\ i=1,\cdots,m.
\end{equation}
Furthermore, we can take $g_i$ to be of the form $p_i\x^{\c}+q_i\y^{\c}$, for $p_i,q_i$ elements of degree $0\in\s$ in $k[\x,\y]$.
\end{lemma}

\begin{proof} The first part follows directly from the inclusion $I(\x,\z)\subset I(\x,\y)+I(\y,\z)$ and the homogeneity of $I$.

For the last part, we can interpret (\ref{coc2}) as an equality in $k[\s^{\o_{\t} 3}]$. More precisely, if we regard each $g_i$ as an element of $k[\s^{\o_{\t}2}]$, then we can rewrite (\ref{coc2}) as
\begin{equation}\label{coc3}
\xi_2(g_i)=\sum_{j=1}^m(a_{ij}\xi_3(g_j)+b_{ij}\xi_1(g_j)),\ i=1,\cdots,m.
\end{equation}
Let $P=r\cdot\x^{\a}\o \y^{\b}$ be a nonzero term in $g_i$ for some $i=1,\cdots,m$ ($r\in k$, $\a,\b\in\s$, $\a+\b=\c$). Then $\xi_2(P)$ is also nonzero since $\xi_2$ is injective, hence one of the monomials occurring on the RHS of (\ref{coc3}) must equal $\x^{\a}\o 1\o \z^{\b}$. If $\x^{\a}\o 1\o \z^{\b}=p\cdot\x^{\a'}\o \y^{\b'}\o 1$ for $p$ a monomial of degree 0 and $\a',\b'\in\s$, then $\a'+\b'=\c$ and
\[\x^{\a}\o \y^{\b}=\mu_{2}^1(\x^{\a}\o 1\o \z^{\b})=\mu_{2}^1(p\cdot\x^{\a'}\o \y^{\b'}\o 1)=\mu_{2}^1(p)\cdot\x^{\c}\o 1,\]
where $\mu_{2}^1(p)$ is a monomial of degree 0. Similarly if $\x^{\a}\o 1\o \z^{\b}=p\cdot 1\o\y^{\a'}\o \z^{\b'}$ then $\x^{\a}\o\y^{\b}=\mu_{3}^2(p)\cdot 1\o\y^{\c}$. It follows that $P$ can be replaced by the product of a term of degree 0 with $\x^{\c}$ or $\y^{\c}$, and therefore we can take the $g_i$'s to be of the form $p_i\x^{\c}+q_i\y^{\c}$ with $p_i,q_i$ elements of degree 0, as desired.

\end{proof}

We now consider the ideals $J(\x,\y)=(\I(\x,\y):\x^{\c})\subset k[\x,\y]$ and $J(\x,\y,\z)=(\I(\x,\y,\z):\x^{\c})\subset k[\x,\y,\z]$. We prove that they are binomial and homogeneous (see also \cite{EiSt}, Cor.1.7). Homogeneity of $J$ is a consequence of the homogeneity of $\I$. To check that $J$ is binomial it suffices to show that $J/\I$ is binomial in the quotient monoid ring $k[\s^{\o_{\t}2}]$ (or $k[\s^{\o_{\t}3}]$). The conclusion now follows from the fact that for any commutative ring $k$ and monoid $M$, and for any element $m\in M$, the ideal $(0:m)\subset k[M]$ is generated by differences $m_1-m_2$ with $mm_1=mm_2$.

We now let $J_0(\x,\y)=J_0(\A)\subset k[T']$ and $J_0(\x,\y,\z)=J_0(\A,\B)\subset k[T'^2]$ be the degree 0 parts of $J(\x,\y)$ and $J(\x,\y,\z)$. They are also binomial, being the degree 0 parts of homogeneous binomial ideals. Notice that $J$ and $J_0$ only depend on the class of $\c$ in $\overline{\s}$. We have

\begin{lemma}\label{lmj0} (a) $J_0(\A)$ is generated by differences $\A^s-1$, hence is the ideal of a closed subgroup scheme $T''$ of $T'$. In particular, $J_0(\x,\y)=J_0(\y,\x)$.

(b) $J_0(\x,\y,\z)$ is invariant under permutations of $\x,\y,\z$.

Since $\I(\x,\y,\z)$ is also invariant under permutations of $\x,\y,\z$, it follows that
\[J_0(\x,\y,\z)=(\I(\x,\y,\z):\x^{\c})_0=(\I(\x,\y,\z):\y^{\c})_0=(\I(\x,\y,\z):\z^{\c})_0.\]
\end{lemma}

\begin{proof} (a) follows from the fact that $J_0$ is binomial and $J_0(\textbf{1})=0$.

To prove part (b) we show that $J_0(\x,\y,\z)=J_0(\x,\y)+J_0(\x,\z)+J_0(\y,\z)$. Once this is proved, the conclusion follows from the symmetry of $J_0$ in (a).

By the explicit description of $J_0$ in part (a), it follows that $J_0(\y,\z)\subset J_0(\x,\y)+J_0(\x,\z)$, so it suffices to prove that $J_0(\x,\y,\z)=J_0(\x,\y)+J_0(\x,\z)$. Since $\I(\x,\y),\I(\x,\z)\subset\I(\x,\y,\z)$ we also get that $J_0(\x,\y),J_0(\x,\z)\subset J_0(\x,\y,\z)$. For the reverse inclusion, observe that
\[J_0(\x,\y,\x)=(\I(\x,\y,\x):\x^{\c})_0=(\I(\x,\y):\x^{\c})_0=J_0(\x,\y)\]
and similarly $J_0(\x,\x,\z)=J_0(\x,\z)$. Also, being binomial and homogeneous, $J_0(\x,\y,\z)$ is generated by binomials of the form $p(\x,\y,\z)=\x^s\y^w\z^{-s-w}-1$. We have
\[
\begin{split}
\x^s\y^w\z^{-s-w}-1&=(\x^{-w}\y^w-1)+(\x^{s+w}\z^{-s-w}-1)\x^{-w}\y^w\\
&=p(\x,\y,\x)+p(\x,\x,\z)\x^{-w}\y^w\in J_0(\x,\y,\x)+J_0(\x,\x,\z)=J_0(\x,\y)+J_0(\x,\z).
\end{split}
\]
This shows that $J_0(\x,\y,\z)\subset J_0(\x,\y)+J_0(\x,\z)$, concluding the proof.

\end{proof}

We will need one more technical lemma for the proof of Theorem \ref{torequiv}:

\begin{lemma}\label{dich} Let $g_i=p_i\x^{\c}+q_i\y^{\c}$, $i=1,\cdots,m$, be as in Lemma \ref{lmcoc2}. Then precisely one of the following holds:

1) $\x^{\c}\equiv\y^{\c}\rm{ mod }\I(\x,\y)$ (modulo changing $\c$ with another representative in its class). In this case $g_i\equiv\y^{\c}e_i\rm{ mod }\I(\x,\y)$, where $e_i$ are elements of degree $0\in\s$ whose images in $k[T'']$ generate the ideal $J_0'$ of a closed subgroup scheme of $T''$.

2) There exists no monomial $p$ of degree $0\in\s$ with the property that $\x^{\c}\equiv p(\x,\y)\y^{\c}\rm{ mod }\I(\x,\y)$. In this case the following congruences hold for $i=1,\cdots,m$:
\begin{equation}\label{relx}
p_i(\x,\z)\equiv\sum_{j=1}^m a_{ij}p_j(\x,\y)\rm{ mod }J_0(\x,\y,\z).
\end{equation}
\begin{equation}\label{rely}
\sum_{j=1}^m a_{ij}q_j(\x,\y)\equiv-\sum_{j=1}^m b_{ij}p_j(\y,\z)\rm{ mod }J_0(\x,\y,\z).
\end{equation}
\begin{equation}\label{relz}
q_i(\x,\z)\equiv\sum_{j=1}^m b_{ij}q_j(\y,\z)\rm{ mod }J_0(\x,\y,\z).
\end{equation}
\end{lemma}

\begin{proof} If 2) doesn't hold, then $\x^{\c}\equiv p(\x,\y)\y^{\c}\rm{ mod }\I(\x,\y)$ for some $p=\x^{s}\y^{-s}$, $s\in l(\s)$. It follows that $\x^{\c-s}\equiv\y^{\c-s}\rm{ mod }\I(\x,\y)$, and since $\c,\c-s$ represent the same class in $\overline{\s}$, 1) must hold. Of course, if 1) holds then 2) doesn't, therefore one and only one of 1) and 2) is true.

Suppose case 1) holds. We can assume that $\x^{\c}\equiv\y^{\c}\rm{ mod }\I$. Then $g_i\equiv\y^{\c}e_i\rm{ mod }\I$, where $e_i=p_i+q_i$. Let $J_0'(\A)$ be the ideal generated by $e_1(\A),\cdots,e_m(\A)$ in $k[T'']=k[T']/J_0$. From (\ref{coc2}) and the fact that $\x^{\c}\equiv\y^{\c}\equiv\z^{\c}\rm{ mod }\I(\x,\y,\z)$, it follows that
\[e_i(\x,\z)\equiv\sum_{j=1}^{m}(a_{ij}e_j(\x,\y)+b_{ij}e_j(\y,\z))\rm{ mod }J_0(\x,\y,\z),\]
which we can rewrite as
\[e_i(\A\B)=\sum_{j=1}^{m}(a_{ij}e_j(\A)+b_{ij}e_j(\B))\rm{ in }k[T''^2],\]
thus $J_0'(\A\B)\subset J_0'(\A)+J_0'(\B)$. Since $e_i(\x,\x)\x^{\c}=g_i(\x,\x)=0$, we get $e_i(\1)=e_i(\x,\x)=0$, whence $J_0'(\1)=0$. The symmetry of $I$ gives the inclusions $g_i(\y,\x)\in (g_j(\x,\y):j=1,\cdots,m)+\I(\x,\y)$, or equivalently $e_i(\A^{-1})\in J_0'(\A)\subset k[T'']$. It follows that $J_0'(\A)=J_0'(\A^{-1})$ and therefore $J_0'$ defines a closed subgroup scheme of $T''$ (see \cite{Wat}, chapter 2), as desired.

Suppose now case 2) holds. We can rewrite (\ref{coc2}) as
\[\x^{\c}\left(p_i(\x,\z)-\sum_{j=1}^m a_{ij}p_j(\x,\y)\right)-\y^{\c}\left(\sum_{j=1}^m a_{ij}q_j(\x,\y)+\sum_{j=1}^m b_{ij}p_j(\y,\z)\right)-\z^{\c}\left(\sum_{j=1}^m b_{ij}q_j(\y,\z)-q_i(\x,\z)\right)\]
\[\in\I(\x,\y,\z).\]
Relations (\ref{relx}-\ref{relz}) then follow as soon as we prove that no two of the monomials $\x^{\c}\o 1\o 1$, $1\o\y^{\c}\o 1$ and $1\o 1\o\z^{\c}$ differ by a monomial of degree 0 in $k[\s^{\o_{\t}3}]$. Suppose $\x^{\c}\o 1\o 1 = q\cdot 1\o\y^{\c}\o 1$ for $q$ a monomial of degree 0, the other cases being analogous. Then
\[\x^{\c}\o 1=\mu_{3}^2(\x^{\c}\o 1\o 1)=\mu_{3}^2(q\cdot 1\o\y^{\c}\o 1)=\mu_{3}^2(q)\cdot 1\o\y^{\c}.\]
Since $q$ has degree 0, the same is true for $\mu_{3}^2(q)$. If we let $p\in k[\x,\y]$ be a monomial of degree 0 representing $\mu_{3}^2(q)$ we get that $\x^{\c}\equiv p(\x,\y)\y^{\c}\rm{ mod }\I(\x,\y)$, a contradiction.

\end{proof}

If we are in the first case of the preceding lemma then $J_0'$ is generated by differences $\x^{s_i}\y^{-s_i}-1$, hence we can assume that $e_i(\x,\y)=\x^{s_i}\y^{-s_i}-1$. It follows that
\[g_i\equiv e_i\y^{\c}\equiv\y^{-s_i}(\x^{s_i}\y^{\c}-\y^{\c+s_i})\equiv\y^{-s_i}(\x^{\c+s_i}-\y^{\c+s_i})\rm{ mod }\I.\]
Since $\y^{-s_i}$ are units, $I$ is generated by the differences $\x^{\c+s_i}-\y^{\c+s_i}$ in $r$-degree at most $\c$ and we are done.

From here on we will assume that we are in the second case of Lemma \ref{dich}. We can rewrite conditions (\ref{relx}-\ref{relz}) as
\begin{equation}\label{rex}
P(\A\B)=A(\A,\B)P(\A)
\end{equation}
\begin{equation}\label{rey}
A(\A,\B)Q(\A)=-B(\A,\B)P(\B)
\end{equation}
\begin{equation}\label{rez}
Q(\A\B)=B(\A,\B)Q(\B)
\end{equation}
where $A=(a_{ij})$, $B=(b_{ij})$, $P=(p_i)$, $Q=(q_i)$, and equality is interpreted as taking place in the coordinate ring of $T''^2$. Since $\x^{\c}$ and $\y^{\c}$ don't differ modulo $\I$ by a monomial of degree 0, it follows that for an element $a(\x,\y)\x^{\c}+b(\x,\y)\y^{\c}$ with $a,b$ elements of degree 0 to be contained in $I$ it is necessary and sufficient that $(a(\A),b(\A))$ is a linear combination of $(p_i(\A),q_i(\A))$ in $k[T'']^2$.

Suppose first that the $p_i$'s don't generate the unit ideal in $k[T'']$. Then $P(u_0)=0$ for some $u_0\in T''$ (after some base change), hence (\ref{rex}) together with the fact that $T''$ is a group shows that $P=0$ in $k[T'']$. The symmetry of $I$ then shows that $q_i(\y,\x)\x^{\c}\in I$, which by the remark in the previous paragraph is the same as $(q_i(\A^{-1}),0)$ being a linear combination of $(0,q_i(\A))$, yielding $q_i=0\in k[T'']$ for all $i$, which is impossible.

We can therefore assume that the $p_i$'s do generate the unit ideal and hence that one of the generators of $I$ in degree $\c$ is $g_1=\x^{\c}+q_1\y^{\c}$. Furthermore, replacing $g_i$ by $g_i-p_ig_1$, we can assume that $p_i=0$ for $i>1$.

Let $\J(\A)$ denote the ideal generated by $q_2(\A),\cdots,q_m(\A)$ in $k[T'']$. Reflexivity of $I$ and the vanishing of the $p_i$'s for $i>1$ imply that $q_i(\1)=0$ for $i>1$, i.e. $\J(\1)=0$. The symmetry of $I$ shows that the pairs $(q_i(\A^{-1}),p_i(\A^{-1}))$ are linear combinations of $(p_i(\A),q_i(\A))$, which translates into
\[
\begin{split}
q_1(\A^{-1})q_1(\A)&\equiv 1\rm{ mod }\J(\A),\\
q_i(\A^{-1})q_1(\A)&\equiv 0\rm{ mod }\J(\A),\rm{ for }i>1.
\end{split}
\]
It follows that $q_1$ is invertible modulo $\J$ and that $\J(\A^{-1})\subset\J(\A)$.

Relation (\ref{rex}) shows that the first column of $A(\A,\B)$ is $P=(1,0,\cdots,0)$. From (\ref{rey}) we get that the first column of $B(\A,\B)$ equals $-q_1(\A)P$ modulo $\J(\A)$, and from (\ref{rez}) that
\[
\begin{split}
q_1(\A\B)&\equiv -q_1(\A)q_1(\B)\rm{ mod }\J(\A)+\J(\B),\\
q_i(\A\B)&\equiv 0\rm{ mod }\J(\A)+\J(\B),\rm{ for }i>1.
\end{split}
\]
This shows that $\J(\A\B)\subset\J(\A)+\J(\B)$, hence $\J$ defines a subgroup scheme $T'''$ of $T''$. It also shows that $-q_1$ is a character of $T'''$. It follows that $\J$ is generated by differences $\A^{s_i}-1$ and $q_1\equiv -\A^{s_1}\rm{ mod }\J$, for some $s_i\in l(\s)$. We can therefore assume that $q_1=-\A^{s_1}$ and $q_i=\A^{s_i}-1$, for $i>1$. We then get that
\[g_1=p_1\x^{\c}+q_1\y^{\c}\equiv\x^{s_1}(\x^{\c-s_1}-\y^{\c-s_1})\rm{ mod }\I\]
and
\[g_i=p_i\x^{\c}+q_i\y^{\c}\equiv\y^{s_1-s_i}(\y^{\c-s_1}\x^{s_i}-\y^{\c-s_1}\y^{s_i})\rm{ mod }\I.\]
This shows that $I$ is generated modulo $\I$ by the differences $\x^{\c-s_1}-\y^{\c-s_1}$, $\x^{\c-s_1+s_i}-\y^{\c-s_1+s_i}$ in $r$-degree $\c$, concluding the proof.

\section{Quotients in the finite toric case}\label{quots}

As a consequence of Theorem \ref{torequiv} we prove the existence of effective geometric quotients for finite toric equivalence relations on affine toric varieties. We will need the following criterion which is due to Olivier:

\begin{lemma}(\cite{Oli})\label{olivier} If $A\to B$ is a ring homomorphism, then the following properties are equivalent:

(i) The sequence $\rm{Spec }B\o_A B\rightrightarrows\rm{Spec }B\to\rm{Spec }A$ is exact in the category of schemes.

(ii) $A$ is the kernel of the map $B\to B\o_A B$, $b\mapsto b\o 1-1\o b$, and the morphism $\rm{Spec }B\to\rm{Spec }A$ is submersive, i.e. it induces the quotient topology on $\rm{Spec }A$.

\end{lemma}

We first give an example showing that an affine toric equivalence relation does not have a categorical quotient in general.

\begin{example}\label{exnonfinite} Consider $X=Y=\bb{A}_k^2$ and the map $X\to Y$ given by $(x_1,x_2)\mapsto(x_1,x_1x_2)$. Let $R=X\times_Y X$, which is a toric equivalence relation, and assume that the categorical quotient $X/R$ exists. Then $X/R=\rm{Spec }A$, where $A=\rm{Ker}\left(\mc{O}_X\overset{\pi_1^*-\pi_2^*}{\longrightarrow}\mc{O}_R\right)$, and $\pi_1,\pi_2$ are the two projections $R\rightrightarrows X$. We have that $A=k[x_1x_2^n,n\geq 0]$, so the map $X\to\rm{Spec }A$ is not submersive because it's not even surjective: the ideal $I=(x_1,x_1x_2-1)\subset A$ is not the unit ideal, but $I\mc{O}_X=\mc{O}_X$.
\end{example}

Nevertheless, the situation is better in the case of a finite toric equivalence relation. We have the following

\begin{corollary}\label{torquot} Let $k$ be a field, $X/k$ an affine toric variety, and $R$ a finite toric equivalence relation on $X$. Then the geometric quotient $X/R$ exists and is effective.
\end{corollary}

\begin{proof} Let $A=\rm{Ker}(\mc{O}_X\to\mc{O}_R)$ and $B=\mc{O}_X$. If $B=k[\s]$ then since the map $\mc{O}_X\to\mc{O}_R$ respects the torus action, it follows that $A=k[\t]$ for some submonoid $\t\subset\s$. If we can prove that $B$ is a finite $A$ module, then $q:X=\rm{Spec }B\to\rm{Spec }A$ must be submersive, and since by Theorem \ref{torequiv} $R$ is induced by $q$, it would follow from Lemma \ref{olivier} that $q$ is a categorical quotient, hence also a geometric quotient.

If we look in $r$-degree zero, we see that $k[l(\s)]\to k[l(\s)]\o_{k[l(\t)]}k[l(\s)]$ is finite, hence $l(\t)$ must have finite index in $l(\s)$. Replacing $\t$ by $\t+l(\s)$ we can therefore assume that $l(\t)=l(\s)$. If we let $C=k[l(\s)]$, $\ss=\s/l(\s)$ and $\tt=\t/l(\t)$, then $\tt\subset\ss$, $A,B$ are $\tt,\ss$-graded $C$-algebras and $B$ is a graded $A$-module.

If we let $\m=\bigoplus_{i\neq 0}A_i$ then $A/\m\simeq C$ and since $B\o_A B$ is finite over $B$ (via the map $b\mapsto b\o 1$), we get by tensoring with $A/\m$ that $B/\m\o_C B/\m$ is a finite $B/\m$-module. But $B/\m$ is a free $C$-module, hence it must have a finite basis over $C$. It follows by the graded version of Nakayama's lemma (using again Lemma \ref{lmords}) that $B$ is a finite $A$-module, concluding the proof.

\end{proof}

%\end{document}
\section*{Acknowledgments} I would like to thank J\'anos Koll\'ar for proposing the original question and David Eisenbud for his guidance throughout the project. I would also like to thank Frank Schreyer for suggesting a simplification of the examples I had of finite noneffective equivalence relations, and Dustin Cartwright and Daniel Erman for helpful conversations.

\end{document}